# GLOBAL OPTIMIZATION OF MULTIVARIABLE FUNCTIONS SATISFYING THE VANDERBEI CONDITION

Natalya Arutyunova*, Aidar Dulliev*, Vladislav Zabotin*

## Abstract

We propose two algorithms for solving global optimization problems on a hyperrectangle with an objective function satisfying the Vanderbei condition (this function is also called an ε-Lipschitz continuous function). The algorithms belong to the class of non-uniform coverings methods. For the algorithms we prove propositions about convergence to an ε-solution in terms of the objective function. We illustrate the performance of the algorithms using several test numerical examples with non-Lipschitz continuous objective functions.

* Department of Applied Mathematics and Informatics, Kazan National Research Technical University named after A.N. Tupolev-KAI, Russian Federation. E-mail: {NKArutyunova@kai.ru, AMDulliev@kai.ru, VIZabotin@kai.ru}.

# ГЛОБАЛЬНАЯ МИНИМИЗАЦИЯ ФУНКЦИЙ МНОГИХ ПЕРЕМЕННЫХ, УДОВЛЕТВОРЯЮЩИХ УСЛОВИЮ ВАНДЕРБЕЯ

Н.К. Арутюнова[*], А.М. Дуллиев[*], В.И. Заботин[*]

В работе предлагаются два алгоритма глобальной минимизации вандербеевой (ε-липшицевой) функции на n-мерном параллелепипеде. Оба они относятся к семейству алгоритмов неравномерных покрытий. Проведено обоснование их сходимости к ε-оптимальному решению. Работа алгоритмов демонстрируется на нескольких тестовых примерах с нелипшицевой целевой функцией. Библ. 10. Табл. 4. Рис. 29.

## 1. Введение

Общеизвестно, что условие Липшица играет важнейшую роль в глобальной оптимизации, обусловленную тем, что с его помощью оцениваются изменения значений функций относительно изменений их аргументов[1]. Существует несколько модификаций этого условия, например, условие Гельдера, полулипшицевости [1], *Q*-липшицевости [2], $\varphi$-$(\alpha,\beta,\nu,\delta,\omega)$-липшицевости [3] и т.д. Одна из таких модификаций в контексте непрерывной одномерной минимизации была предложена в 1999 году Вандербеем в статье [4]. (Ниже она называется условием Вандербея, а функция, ему удовлетворяющая, – непрерывной по Вандербею, или короче, вандербеевой[2].) В последующие годы условие Вандербея также использовалось при построении численных методов решения различных оптимизационных задач, например, при обобщении метода Ю.Г. Евтушенко [5] для задачи минимизации вандербеевой функции на отрезке [6], в задаче поиска проекции заданной точки на множество решений уравнения $f(x)=0$ [7], в задаче глобальной липшицевой оптимизации при наличии шума [8].

В данной работе предлагаются два алгоритма глобальной минимизации вандербеевой функции на *n*-мерном параллелепипеде. Оба они относятся к семейству методов неравномерных покрытий и отличаются от других его представителей (см., например, [5]), во-первых, способом построения элементов покрытия, базирующемся на условие Вандербея, и, во-вторых, способом размещения элементов покрытия в параллелепипеде.

---

[*] Кафедра прикладной математики и информатики, Казанский национально-исследовательский технический университет им. А.Н. Туполева, Российская Федерация. E-mail: {NKArutyunova@kai.ru, AMDulliev@kai.ru, VIZabotin@kai.ru}.

[1] Исследованию задач липшицевой глобальной оптимизации посвящена обширная литература, подробную информацию о которой можно найти, например, в [8, 9]

[2] В литературе можно встретить другое название функции, удовлетворяющей условию Вандербея,– $\varepsilon$-липшицева функция, а также – слегка видоизмененное понятие – $\delta$-липшицева функция. В настоящей работе авторы не придерживаются этих названий, так как символы $\varepsilon$ и $\delta$ используются для записи точностей решения рассматриваемой задачи.

Общая идея первого метода заключается в последовательном рекурсивном разбиении рассматриваемого параллелепипеда в зависимости от его размеров на ($n$ +1) или менее новых параллелепипедов (вдоль осей пространства) и сравнении значений функции в некоторых их точках. Количество параллелепипедов разбиения, а следовательно, и быстродействие метода существенно зависит от порядка перебора параллелепипедов. Соответствующие схемы перебора наглядно описываются с помощью дерева параллелепипедов разбиения, которое при программной реализации просто описывается с помощью списков.

Второй метод является разновидностью метода ветвей и границ. На каждой итерации либо производится исключение из дальнейшего рассмотрения текущего параллелепипеда («отсечение»), либо текущий параллелепипед разбивается пополам по ребру максимальной длины («ветвление»), либо из текущего параллелепипеда удаляется некоторый параллелепипед, центр которого расположен в центре текущего, и затем, оставшееся множество дробится по определенной схеме на $2n$ параллелепипедов («отсечение плюс ветвление»). Размеры удаляемого параллелепипеда вычисляются с помощью нижней оценки значений целевой функции в текущем параллелепипеде, получаемой из условия Вандербея с использованием текущей оценки глобального минимума (рекорда). Выбор способа дробления осуществляется в зависимости от размеров удаляемого параллелепипеда в схеме «отсечение плюс ветвление», причем таким образом, чтобы по возможности избежать появления параллелепипедов, имеющих слишком малую длину какого-либо ребра.

Предлагаемые алгоритмы тестируются на нескольких примерах с нелипшицевой целевой функцией.

## 2. Постановка задачи и описание алгоритмов

Пусть на множестве $A \subset \mathbb{R}^n$ задана функция $f : A \to \mathbb{R}$ и $\| \cdot \|$ – евклидова норма в $\mathbb{R}^n$.

**Определение 1** [4]. *Функция $f(x)$, называется непрерывной по Вандербею на A, если для любого $\eta > 0$ существует $L(\eta) > 0$, что для всех $x, y \in A$ выполняется условие*

$$| f(x) - f(y) | \le L(\eta) \| x - y \| + \eta. \tag{1}$$

Некоторые свойства функции $L(\eta)$, а также способы ее нахождения, можно найти в работах [4, 7, 10].

Рассмотрим задачу

$$f_* := \operatorname*{glob\,min}_{x \in P} f(x) \tag{2}$$

где $P = \{x \in \mathbb{R}^n \mid a^i \le x^i \le b^i, i = 1,...,n\}$ – $n$-мерный параллелепипед, $f(x)$ – непрерывна по Вандербею на $P$.

Приведем два алгоритма поиска ε-оптимального решения задачи (2), то есть поиска какой-либо точки из множества

$$X_*^{\varepsilon} = \left\{ x \in P \mid f(x) \le f_* + \varepsilon \right\}. \quad (3)$$

Оба алгоритма относятся к семейству алгоритмов неравномерных покрытий и отличаются от других его представителей, во-первых, способом построения элементов покрытия, базирующемся на условие (1); во-вторых, способом размещения элементов покрытия в параллелепипеде $P$.

В алгоритмах 1 и 2 будут использоваться следующие общие обозначения: $k$ – номер итерации; $\tilde{x}_k := \arg\min\left\{ f(x) \mid x \in \{x_1, ..., x_k\} \right\}$, $F_k := f(\tilde{x}_k)$; $\mathcal{P}_k$ – текущий набор параллелепипедов $P_i$, генерируемых алгоритмом.

Все параллелепипеды считаются замкнутыми.

**Алгоритм I**

**Шаг 0.** Задать $a_0 = a = (a^1, \ldots, a^n)$ и $b_0 = b = (b^1, \ldots, b^n)$ – вершины рассматриваемого параллелепипеда $P_0 = X$, а также величины $\varepsilon > 0$ и $\eta \in (0; \varepsilon)$.

Принять $F_{-1} = f(a)$; $\mathcal{P}_0 = \{P_0\}$, $k = 0$, $m = 1$.

Вычислить:

$L(\eta)$ – по известной априорно зависимости или с помощью одного из алгоритмов оценки;

$h$ – базовую величину шага по формуле:

$$h = \frac{2(\varepsilon - \eta)}{L(\eta)}. \quad (4)$$

**Основной цикл.** Повторять, пока набор $\mathcal{P}_k$ не пуст.

Работа с одним параллелепипедом $P_k$ из набора $\mathcal{P}_k$.

**Шаг 1.** Взять новую итерационную точку $x_k$:

$$x_k^i = \min\left\{ a_k^i + \frac{h}{2}; b_k^i \right\}, \; i = 1, ..., n,$$

и значение $f(x_k)$.

**Шаг 2.** Обновить при необходимости известное на данный момент минимальное значение функции $f$ или скорректировать шаг:

*если* $f(x_k) > F_{k-1}$, *то*

- $F_k = F_{k-1}$;
- вычислить актуальный шаг:

$$h' = h + \frac{f(x_k) - F_k}{L(\eta)};$$

*иначе* $F_k = f(x_k)$ и $h' = h$.

**Шаг 3.** Произвести разбиение $P_k$ на серию параллелепипедов $P_{k,i}$ ($i$ – номер координаты, по которой производится разбиение, $i \le n$).

*Цикл по* $i = \overline{1, n}$:

если $b_k^i - a_k^i > h'$, ● ввести новый параллелепипед $P_{k,i}$ с углами:

$$a_{k,i} = \left(a_k^1, \ldots, a_k^{i-1}, a_k^i + h', a_k^{i+1}, \ldots, a_k^n\right);$$

$$b_{k,i} = \left(\min\left\{a_k^1 + h', b_k^1\right\}, \ldots, \min\left\{a_k^{i-1} + h', b_k^{i-1}\right\}, b_k^i, \ldots, b_k^n\right);$$

● $m := m + 1$.

(На параллелепипеде $P_{k,n+1}$ с углами: $a_{k,n+1} = a$ и $b_{k,n+1}$: $\forall i = \overline{1, n}\left(b_k^i = \min\left\{a_k^i + h'; b_k^i\right\}\right)$, ближайшем к началу координат, для функции $f$ гарантированно выполнится условие:

$$\forall x \in P_{k,n+1} : f(x) \le F_k - \varepsilon).$$

**Шаг 4.** Исключить из набора рассмотренный параллелепипед $P_k$ и добавить все созданные на шаге 3 параллелепипеды $P_{k,i}$, $i \le n$. Полученный список назовём $\mathcal{P}_{k+1}$.

**Шаг 5.** Принять: $m := m - 1$, $k := k + 1$ и перейти к новой итерации основного цикла.

Алгоритм завершает свою работу, когда список $\mathcal{P}_k$ оказывается пуст, при этом точка $\tilde{x}_{k-1}$ принимается в качестве ε-оптимального решения задачи (2).

***Комментарии к алгоритму I***

**1.** При программной реализации алгоритма наборы параллелепипедов можно хранить в виде списков. Самыми удобными и, обычно, быстрыми операциями по выбору, удалению и добавлению элементов являются действия с началом (головой) и концом списка. Комбинация действий (*с началом* или *концом списка*) при выборе параллелепипеда для разбиения (а, значит, при его *удалении из списка*) и при *добавлении в список* новых параллелепипедов будет определять порядок расположения параллелепипедов в списке и порядок их перебора. Поскольку размеры новых параллелепипедов разбиения, а значит, и скорость процесса покрытия всего исходного параллелепипеда, существенно зависят от известного на данный момент лучшего среди всех рассмотренных итерационных точек значения функции (шаг 2), порядок перебора будет напрямую влиять на объёмы оперируемых данных и длительность вычислений.

Случай, когда следующий для разбиения параллелепипед выбирается из набора самых последних параллелепипедов, добавленных в список (на предыдущей итерации основного цикла), соответствует обходу дерева параллелепипедов разбиения *в глубину* (схема 1). При выборе же параллелепипеда для разбиения из группы, добавленной в список раньше всех остальных (на предыдущем уровне разбиения), получается обход дерева параллелепипедов разбиения *в ширину* (схема 2).

Получить эти схемы можно следующими основными комбинациями операций с началом и концом списка (см. табл. 1).

| № варианта | Действие | | | Схема |
|---|---|---|---|---|
| | выбор (удаление) $P_k$ из $\mathcal{P}_k$ | добавление $P_{k,i}$ в новый подсписок ($\mathcal{P}'_k$) | добавление $\mathcal{P}'_k$ в $\mathcal{P}_k$ | |
| 1 | из начала | в начало | в начало | 1а |
| 2 | из начала | в конец | в начало | 1б |
| 3 | из начала | в начало | в конец | 2а |
| 4 | из начала | в конец | в конец | 2б |

Таблица 1: Основные варианты реализации алгоритма 1

Для каждой из приведённых схем можно сформировать симметричный вариант – он будет отличаться лишь фактическим расположением элементов в списке, но не самим принципом разбиения, поэтому эквивалентен указанной основной схеме.

На рис. 1 и 2 приведены иллюстрации для схем 1 и 2 обхода дерева параллелепипедов разбиения и используемых для этого списков на примере вариантов 1б и 2б (см. табл. 1).

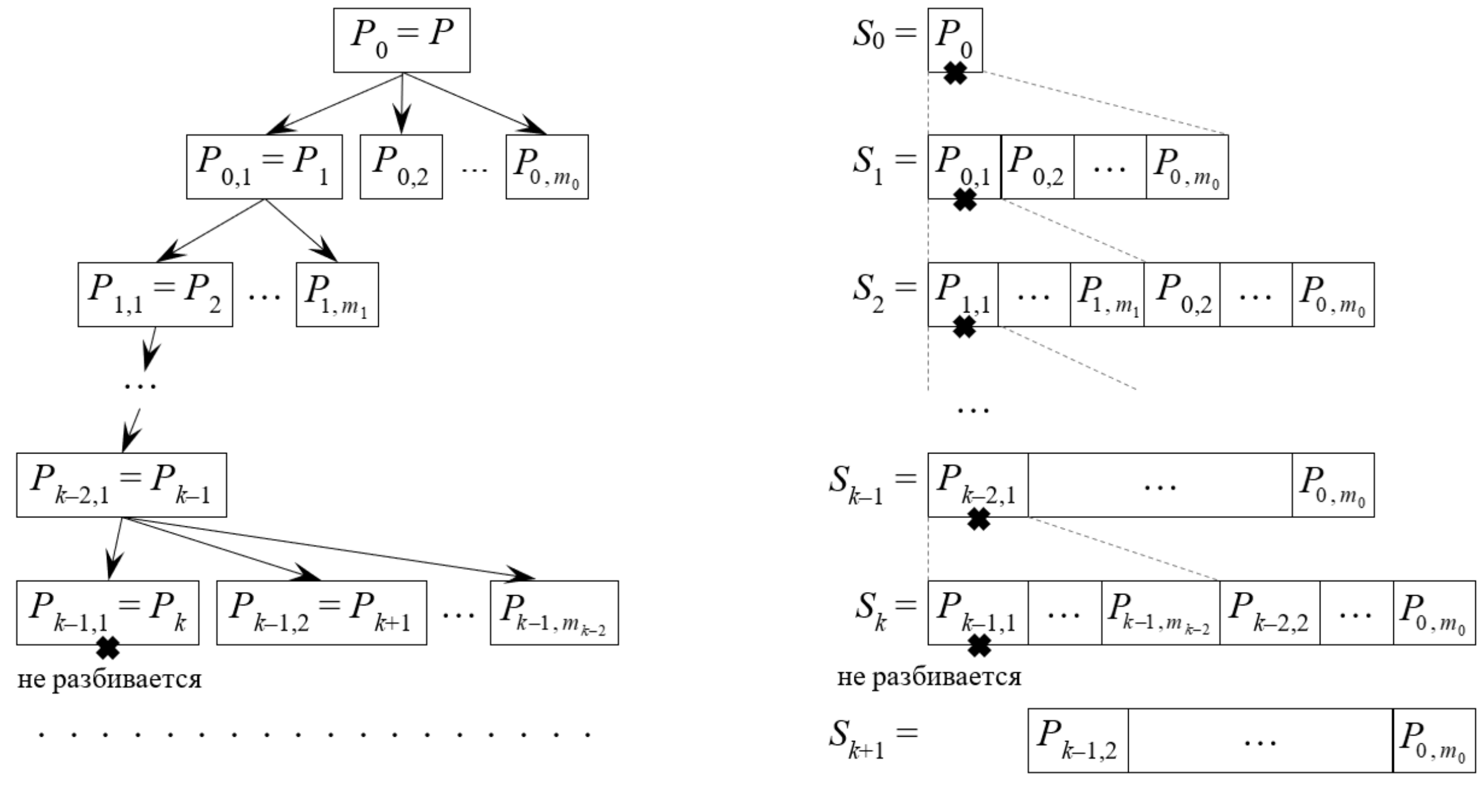

а) дерево параллелепипедов разбиения б) списки параллелепипедов разбиения

Рис. 1. Схема 1б: обход дерева параллелепипедов разбиения в глубину

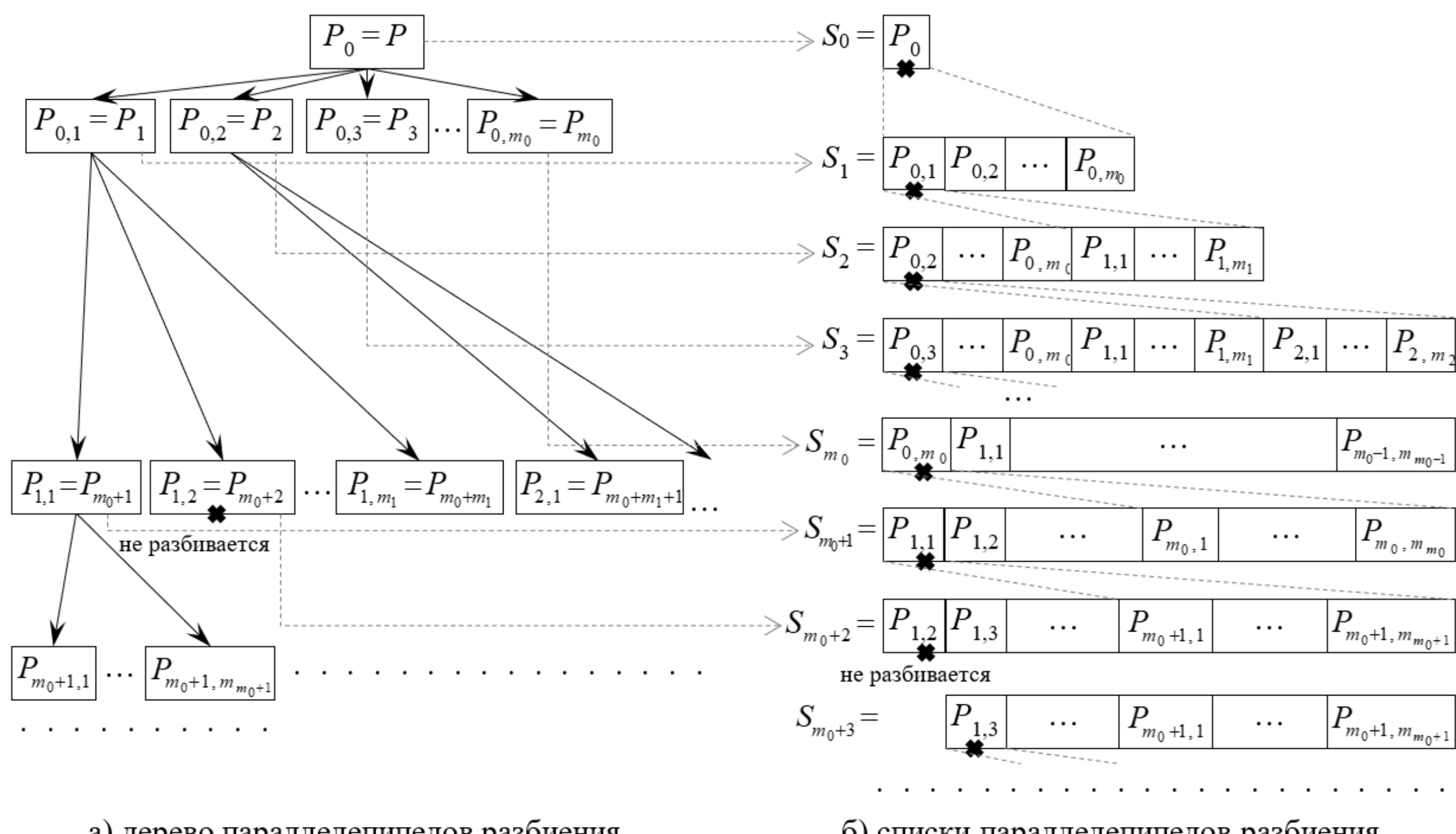

Рис. 2. Схема 2: обход дерева параллелепипедов разбиения в ширину

Для приведённых в табл. 1 четырёх схем обхода были проведены численные эксперименты, их результаты и анализ полученных данных приведены ниже в разделе 3.

**2.** Описанный алгоритм допускает возможность распараллелить процесс поиска решения – рассмотрение различных параллелепипедов набора $\mathcal{P}_k$ одновременно (несколькими процессорами), добавление затем образованных подсписков в общий набор, определение «лучшего» значения функции (общего) и использование его для новых параллельных расчётов.

При последовательной же реализации этого алгоритма есть возможность использовать рекурсивный поиск: при создании нового параллелепипеда $P_{k,i}$ сразу его рассматривать (и при необходимости разбивать). Такой порядок рассмотрения, очевидно, соответствует описанной схеме 1б. Преимуществом такого варианта реализации является отсутствие необходимости использования списков: работа всегда ведётся с последним созданным, но не рассмотренным параллелепипедом. При этом не только не требуется выделение дополнительной памяти для хранения списка, но также сокращаются затраты времени: не нужно выполнять дорогостоящие операции со списками. Однако нужно помнить о том, что при переходе к новому шагу рекурсии все прежние данные сохранятся в стек, размер которого ограничен, поэтому возможность использования данного подхода зависит от сложности решаемой задачи.

В приведённых ниже результатах расчётов для схемы 1б для примеров, где оказался возможен поиск решения рекурсивным способом, указывается именно соответствующее ему время расчётов. В случае же невозможности решения с применением рекурсии (из-за переполнения стека), расчёты производились с использованием списков.

Сравнение же этих двух вариантов реализации схемы 1б показало, что рекурсивный подход в случае применимости даёт сокращение затрат времени примерно в 10 раз по сравнению со списковым подходом.

**3.** Величины шага $h$ и $h'$ существенно зависят от параметра η. При этом, судя по (4), для получения положительной величины шага требуется выполнение соотношения: $\eta < \varepsilon$. Удачный выбор соотношения величин η и ε даёт значительное снижение затрат на расчёты, но для каждой решаемой задачи он индивидуален.

Поэтому при проведении расчётов для выбранной задачи были сначала произведены эксперименты с разными вариантами значений η и ε для определения оптимального их соотношения (с использованием наиболее быстрого, рекурсивного подхода), а затем именно оно использовалось при выполнении следующих тестов.

### Алгоритм II

Приводимый ниже алгоритм излагается для случая, когда $\inf\{L(\eta)\mid \eta>0\}>0$.

Дополнительно к введённым ранее обозначениям при описании этого алгоритма также используются следующие: $c(P_i)$, $d_i$ – центр и длина диагонали параллелепипеда $P_i$ соответственно; γ – параметр метода, $r$ – половина длины диагонали параллелепипеда $P$. Условие, при котором выполняется соответствующая операция, записывается в квадратных скобках.

**Шаг 0.** Положить $P_1 = P$, $x_1 = c(P_1)$, $\mathcal{P}_1 = \{P_1\}$, $\mu_1 = 0$. Задать параметры метода: $\varepsilon > 0$, $\beta \in (0; 1)$, $\gamma \in (r_1 / r; 1]$, вычислив $r_1$ по формуле (5). Установить $k = 1$.

**Шаг 1.** Вычислить

$$r_k = \min\{\rho_k; r\}, \tag{5}$$

где $\rho_k$ находится как решение задачи

$$\rho_k = \sup\left\{\frac{f(x_k)-F_k+\varepsilon-\eta}{L(\eta)} \mid 0<\eta\le\varphi_k := f(x_k)-F_k+\beta\varepsilon\right\}.$$

**Шаг 2.** Последовательно проверять условия, приведенные ниже в квадратных скобках, и при первом выполнении одного из них осуществить соответствующую операцию, после чего, перейти на шаг 3.

*2.1.* $[r_k \ge d_k/2]$. Положить $\mu_{k+1} = \mu_k + \mathrm{vol}(P_k)$. Построить новый набор $\mathcal{P}_{k+1}$ из набора $\mathcal{P}_k$, удалив в последнем $P_k$.

*2.2.* $\left[r_k < d_k/2 \wedge r_k < \gamma r\right]$. Параллелепипед $P_k$ разделить пополам гиперплоскостью, перпендикулярной какому-либо максимальному ребру. Построить новый набор $\mathcal{P}_{k+1}$ из набора $\mathcal{P}_k$, добавив к нему полученные параллелепипеды и удалив $P_k$. Положить $\mu_{k+1} = \mu_k$ .

*2.3.* $\left[r_k < d_k/2 \wedge r_k \geq \gamma r\right]$. Построить параллелепипед $\hat{P}_k$, удовлетворяющий условиям:

а) $\hat{P}_k \subset P_k \cap B$, где $B = \left\{x \in \mathbb{R}^n \mid \| x - x_k \| \leq r_k\right\}$ ;

б) грани $\hat{P}_k$ параллельны координатным плоскостям;

в) $\hat{P}_k$ имеет максимальный объем $\mu$.

Положить $\mu_{k+1} = \mu_k + \mathrm{vol}(\hat{P}_k)$. Построить новый набор $\mathcal{P}_{k+1}$, удалив из набора $\mathcal{P}_k$ параллелепипед $P_k$ и добавив параллелепипеды, получаемые следующим образом.

А. Разбить $P_k$ на три параллелепипеда, проведя две гиперплоскости через грани параллелепипеда $\hat{P}_k$, перпендикулярные какому-либо максимальному ребру параллелепипеда $P_k$.

Б. Из трех полученных в $P_k$ параллелепипедов выбрать два, не содержащие $\hat{P}_k$, и добавить их в $\mathcal{P}_{k+1}$.

В. Провести операции, аналогичные А, Б, заменив $P_k$ параллелепипедом из п. А, содержащим $\hat{P}_k$.

Г. Повторять операции А–В, до тех пор, пока параллелепипед из п. А, содержащий $\hat{P}_k$, не совпадет с $\hat{P}_k$.

**Шаг 3.** Среди всех параллелепипедов, входящих в $\mathcal{P}_{k+1}$, выбрать

$$P_{k+1} := \arg\min \{ f(c(P_i)) \mid P_i \in \mathcal{P}_{k+1} \}$$

и положить $x_{k+1} = c(P_{k+1})$.

**Шаг 4.** Проверить условие останова: если $\mu_{k+1} = \mathrm{vol}(P)$, то принять точку $\tilde{x}_{k+1}$ в качестве $\varepsilon$-оптимального решения задачи (2) и завершить работу алгоритма. Иначе перейти к шагу 1, приняв $k := k + 1$.

***Комментарии к алгоритму II***

**1.** Очевидно, на каждой $k$–й итерации $\rho_k \geq \dfrac{f(x_k) - F_k + (1-\beta)\varepsilon}{L(\beta\varepsilon)} \geq \dfrac{(1-\beta)\varepsilon}{L(\beta\varepsilon)} > 0$; значит, $r_k \geq \min\left\{\dfrac{(1-\beta)\varepsilon}{L(\beta\varepsilon)}, r\right\}$. Отсюда и из неравенства $r_k \leq r$ вытекает, что величина $r_k$ определяется корректно. Точное нахождение $r_k$ возможно далеко не всегда, поэтому на практике можно ограничиться приближенным вычислением. При этом на каждой $k$–й итерации целесообразно пользоваться информацией о значениях $r_i$, $\varphi_i$, $i = 1, ..., k - 1$, так как, во-первых, среди них может найтись $\varphi_j$, совпадающее с $\varphi_k$, и тогда $r_k = r_j$, во-вторых, наибольшее значение среди $\varphi_i$, $i = 1, ..., k - 1$, позволяет сузить область поиска $r_k$.

**2.** Параметр $\gamma$, используемый на операциях *2.2*, *2.3*, влияет на выбор алгоритмом схемы разбиения текущего параллелепипеда $P_k$. Если $r_k < \gamma r$, то осуществляется половинное деление $P_k$; иначе из $P_k$ удаляется параллелепипед $\hat{P}_k$, на котором достигнутый рекорд $F_k$ не может быть улучшен более чем на $\varepsilon$, а затем оставшееся множество разбивается на $2n$ параллелепипедов, каждый из которых имеет ровно одну грань, содержащую одну из $2n$ граней параллелепипеда $\hat{P}_k$. Поскольку значение $\gamma$ определяет нижнюю границу диагонали удаляемого параллелепипеда $\hat{P}_k$, желательно, чтобы оно не было настолько малым, что после удаления $\hat{P}_k$ и дальнейшего разбиения $P_k \setminus \hat{P}_k$ какие-то ребра получающихся параллелепипедов были бы очень «узкими» ($\approx \delta := 2r_1 / \sqrt{n}$), а также, чтобы оно не было настолько большим, что удаление $\hat{P}_k$ происходило бы довольно редко (очевидно, что при $\gamma = 1$ производится только половинное деление параллелепипедов $P_k$). По-видимому, выбор конкретных значений $\gamma$, обеспечивающих наилучшую скорость сходимости, весьма трудная задача. Анализ экспериментальных данных, полученных авторами на ряде численных примеров при различных $\gamma \in (r_1 / r;\ 1]$, показал, что количество вычисленных значений целевой функции для некоторых $\gamma < 0.1$ может оказаться примерно на 10–40% меньше, чем для $\gamma = 1$ (см. раздел 3).

**3.** На шаге 3 поиск минимума можно осуществлять следующим образом. Список $\mathcal{P}_k$ изменяется на каждой итерации таким образом, что либо только удаляется параллелепипед (операция *2.1*), либо вначале удаляется параллелепипед, а затем, добавляются два (операция *2.2*) или $2n$ (операция *2.3*) параллелепипедов. Этот набор целесообразно хранить в виде упорядоченного по возрастанию $f(c(P_i))$ списка, в котором на каждой итерации добавляются (удаляются) параллелепипеды в соответствие с данным порядком.

**Предложение 1.** *Условие останова* $\mu_{k+1} = \mathrm{vol}(P)$ *на шаге 4 достигается за конечное число итераций.*

**Доказательство.** Изменение величины $\mu_k$ происходит либо на операции *2.1*, либо на *2.3*. На каждой из них она увеличивается на объем параллелепипеда, точки которого, не совпадающие с его центральной точкой, исключаются из исследования на всех последую-щих итерациях. Все эти параллелепипеды содержаться в $P$, попарно не имеют общих внутренних точек и покрывают $P$, следовательно, $\sum_{k\geq 1}(\mu_{k+1}-\mu_k)=\mathrm{vol}(P)$. Для доказательства предложения достаточно убедиться в справедливости неравенства $\inf_{k\geq 1}(\mu_{k+1}-\mu_k)>0$.

Пусть $\mathcal{Q}$ – множество всех параллелепипедов, получаемых во всех итерациях на операциях *2.1–2.3* и при этом входящих в наборы $\mathcal{P}_k$ или совпадающих с $\hat{P}_k$, т.е. $\mathcal{Q}=\left\{\hat{P}_k \mid k\geq 1, \hat{P}_k \textit{ существует}\right\}\cup\bigcup_{k\geq 1}\mathcal{P}_k$. Упорядочим ребра каждого $P_j\in\mathcal{Q}$ по возрастанию: $a_j^{(1)}\leq a_j^{(2)}\leq\ldots\leq a_j^{(n)}$.

Рассмотрим следующие множества

$\mathcal{Q}^0=\left\{P_j\in\mathcal{Q}\mid a_j^{(1)}>\delta\right\}$, $\mathcal{Q}^i=\left\{P_j\in\mathcal{Q}\mid a_j^{(i)}\leq\delta, a_j^{(i+1)}>\delta\right\}$, $i=1,\ldots,n-1$, $\mathcal{Q}^n=\left\{P_j\in\mathcal{Q}\mid a_j^{(n)}\leq\delta\right\}$.

Покажем, что все они конечны.

Конечность $\mathcal{Q}^0$ очевидна. По построению и в силу неравенства $\gamma>r_1/r$ для каждого $i=1,\ldots,n$ и для любого $P_\lambda\in\mathcal{Q}^i$ найдется $P_\nu\in\mathcal{Q}^{i-1}$ такой, что $P_\lambda$ получается из $P_\nu$ в результате конечного числа дроблений последнего. Параллелепипед $P_\nu$ содержит конечное число параллелепипедов из $\mathcal{Q}^i$. Действительно, $P_\lambda$ получается из $P_\nu$ в результате дробления по некоторому ребру $a_\nu^{(q)}, q\geq i$, на одной какой-нибудь итерации и возможных дроблений на последующих за ней итерациях по каким-либо другим ребрам $a_\nu^{(q')}, q'\geq i, q'\neq q$, которые в свою очередь после дробления не могут оказаться меньше $\delta$, так как $\gamma>r_1/r$. В частности, отсюда следует, что количество параллелепипедов из $\mathcal{Q}^i$, содержащихся в $P_\nu$, будет не более $(2n+1)\prod_{s=i}^{n}\left(\left[\frac{a_\nu^{(s)}}{\delta}\right]+1\right)$. Таким образом, по индукции получаем конечность каждого множества $\mathcal{Q}^i$, $i=1,\ldots,n$.

Конечность множеств $\mathcal{Q}^i$ и равенство $\bigcup_{i=0}^{n}\mathcal{Q}^i=\mathcal{Q}$ гарантируют существование $a_*=\min\left\{a_j^{(i)}\mid P_j\in\mathcal{Q}, i=1,\ldots,n\right\}$. Значит, на операциях *2.1* и *2.3* получаем соответственно $\mathrm{vol}(P_k)\geq a_*^n$ и $\mathrm{vol}(\hat{P}_k)\geq a_*^n$, т.е. $\mu_{k+1}-\mu_k\geq a_*^n>0$ при любых $k\geq 1$. Последнее, как указывалось выше, доказывает предложение. □

Пусть $k_* < \infty$ – номер итерации, в которой на шаге 4 согласно доказанному предложению 1 выполнятся условие останова $\mu_{k_*+1} = \mathrm{vol}(P)$.

**Предложение 2.** *Точка* $\tilde{x}_{k_*+1}$, *найденная на шаге 4, принадлежит множеству* $X_*^{\varepsilon}$ *(см. (3)).*

**Доказательство.** Возьмем произвольную точку $x \in P$. Как указывалось выше (см. доказательство предложения 1), семейство всех параллелепипедов $P_k$ и $\hat{P}_k$, исключаемых из дальнейшего рассмотрения на операциях *2.1* и *2.3* соответственно, покрывает $P$, следовательно, найдется такой номер итерации $j < \infty$, что эта точка содержится в одном из этих параллелепипедов, $P_j$ или $\hat{P}_j$. Не ограничивая общности, положим $x \in P_j$.

Очевидно, что $\| x - x_j \| \le \rho_j$. Из условия Вандербея для любого $\eta \in (0, \varphi_j]$ имеем

$$f(x_j) \le f(x) + L(\eta) \| x - x_j \| + \eta \le f(x) + L(\eta)\rho_j + \eta \le f(x) + f(x_j) - F_j + \varepsilon .$$

Отсюда $F_j \le f(x) + \varepsilon$. Принимая во внимание произвольность выбора точки $x$ и неравенство $F_{k_*+1} \le F_j$, приходим к оценке $F_{k_*+1} \le f_* + \varepsilon$, т.е. $\tilde{x}_{k_*+1} \in X_*^{\varepsilon}$. □

### 3. Численные эксперименты

Тестирование предложенных алгоритмов было проведено на четырех примерах в $\mathbb{R}^2$ (см. табл. 2). Целевые функции во всех примерах непрерывны по Вандербею, но нелипшицевые на $P$, причем во втором и третьем примерах – это нелипшицевые модификации функции Экли и табличной функции Хольдера соответственно. Оценки вандербеевой функции $L(\eta)$ по норме $\| \cdot \|_1$ приведены в табл. 2. Реализация алгоритмов осуществлялась по евклидовой норме.

| № | Целевая функция | Оценка $L(\eta)$ | Параллелепипед |
|---|---|---|---|
| 1 | $f_1(x, y) = -10 e^{-\sqrt{0.5(\lvert x\rvert + \lvert y\rvert)}}$ | $L(\eta) = \dfrac{25}{2\eta}$ | $[-2,12]^2$ |
| 2 | $f_2(x, y) = f_1(x, y) - e^{0.5(\cos(2\pi x) + \cos(2\pi y))}$ | $L(\eta) = \dfrac{25}{2\eta} + \pi e$ | $[-2,12]^2$ |
| 3 | $f_3(x, y) = -\left\lvert \cos x \cos y \, e^{0.5\left\lvert 1 - \sqrt{\lvert x\rvert + \lvert y\rvert}\right\rvert} \right\rvert$ | $L(\eta) = e^{\alpha/2} + \dfrac{e^{\alpha}}{16\eta}$, $\alpha = \sqrt{20} - 1$ | $[-10,10]^2$ |
| 4 | $f_4(x, y) = \sin(5y) \arcsin x - \sin(5x) \arcsin y$ | см. (6) | $[-1,1]^2$ |

Таблица 2: Тестовые примеры

$$L(\eta) = 5\pi + \begin{cases} \dfrac{2}{\sqrt{1-\tau^2(\eta/2)}}, & \text{если } 0 < \eta/2 < \tilde{\eta}, \\ \pi - \eta/2, & \text{если } \tilde{\eta} \le \eta/2 < \pi, \end{cases} \qquad (6)$$

где $\tilde{\eta} = \pi/2 - \sqrt{\dfrac{1-\sigma}{1+\sigma}} - \arcsin\sigma$; $\sigma$ – корень уравнения $(\pi/2 + \arcsin\sigma)\sqrt{1-\sigma^2} = 1+\sigma$ на промежутке [0; 1); $\tau(\alpha)$ при фиксированном $\alpha$ – корень уравнения $(\pi/2 - \alpha - \arcsin\tau)\sqrt{1-\tau^2} = 1-\tau$ на промежутке [0; 1). (Нетрудно показать, что указанные корни на промежутке [0; 1) – единственные).

Всюду далее использованы следующие обозначения: $N_{tot}$ – общее число параллелепипедов разбиения; $N_{opt}$ – номер параллелепипеда, где был достигнут оптимум.

Как было отмечено, в комментарии 3 к алгоритму I, расчёты проводились в два этапа. Результаты первого из них, направленного на поиск лучшего соотношения величин $\eta$ и $\varepsilon$, проиллюстрированы графиками (рис. А.1, А.2)[3] зависимости числа итераций от значения $\eta$ при фиксированном $\varepsilon = 0{,}5$.

Из рис. А.1, А.2 видно, что для рассматриваемых функций лучше использовать следующие соотношения:

для $f_1$: $\eta = 0{,}9\varepsilon$; для $f_2$: $\eta = 0{,}8\varepsilon$; для $f_3$: $\eta = 0{,}6\varepsilon$; для $f_4$: $\eta = 0{,}5\varepsilon$;

Именно с этими соотношениями проводились расчёты на втором этапе, поэтому всюду далее для краткости из этих двух параметров указывается только значение $\varepsilon$.

В табл. 3 приводятся результаты численных экспериментов при различных $\varepsilon$ для разных схем обхода дерева разбиений. Прочерк означает отсутствие результатов по причине больших затрат времени и /или памяти, значительно превышающих таковые при использовании других схем.

[3] Здесь и далее приводятся ссылки на рисунки, размещённые в Приложении А

| $i$ | ε | Схема обхода | $(x_\varepsilon, y_\varepsilon)$ | $f_i(x_\varepsilon, y_\varepsilon)$ | $N_{tot}$ | $N_{opt}$ | Рис. А.№ |
|---|---|---|---|---|---|---|---|
| 1 | 0.5 | 1а | (–0.0011; 0.0008) | –9.6944 | 603993 | 588784 | 3 |
| | | 1б | (–0.0002; –0.0008) | –9.7829 | 1156717 | 1127196 | |
| | | 2а | (–0.0002; –0.0002) | –9.8596 | 105214288 | 104702038 | 4 |
| | | 2б | (–0.0002; –0.0002) | –9.8596 | 102526635 | 102086326 | |
| | 0.1 | 1а | (0.0000; 0.0001) | –9.9363 | 102764377 | 101177237 | 5 |
| | | 1б | (0.0001; 0.0001) | –9.9130 | 226120051 | 222490998 | |
| | | 2а | – | – | – | – | – |
| | | 2б | – | – | – | – | |
| 2 | 0.5 | 1а | (0.0000; 0.0015) | –12.4467 | 121876 | 111959 | 6 |
| | | 1б | (0.0003; –0.0001) | –12.5738 | 201996 | 192015 | |
| | | 2а | (0.0000; 0.0020) | –12.4072 | 398611 | 355408 | 7 |
| | | 2б | (–0.0011; –0.0013) | –12.3823 | 384541 | 348241 | |
| | 0.1 | 1а | (0.0001; –0.0000) | –12.6490 | 20440621 | 20054425 | 8 |
| | | 1б | (–0.0000; –0.0001) | –12.6491 | 21352428 | 20784419 | |
| | | 2а | (0.0000; –0.0001) | –12.6415 | 445121567 | 435245855 | 9 |
| | | 2б | (0.0000; –0.0000) | –12.6705 | 366830725 | 359820100 | |
| 3 | 0.5 | 1а | (–9.4786; 9.4868) | –5.3339 | 36503 | 7039 | 10 |
| | | 1б | (9.4820; –9.4842) | –5.3340 | 18602 | 8690 | |
| | | 2а | (9.4796; 9.4874) | –5.3339 | 40141 | 39670 | 11 |
| | | 2б | (–9.4830; –9.4865) | –5.3340 | 35567 | 1137 | |
| | 0.1 | 1а | (–9.4821; 9.4824) | –5.3340 | 4424905 | 1810484 | 12 |
| | | 1б | (–9.4829; 9.4826) | –5.3340 | 3983228 | 3857481 | |
| | | 2а | (–9.4821; –9.4819) | –5.3340 | 8164034 | 2349617 | 13 |
| | | 2б | (9.4821; 9.4824) | –5.3340 | 6451383 | 6367202 | |
| 4 | 0.5 | 1а | (0.3402; 1.0000) | –1.8904 | 471 | 342 | 14 |
| | | 1б | (0.3451; 1.0000) | –1.8899 | 446 | 180 | |
| | | 2а | (1.0000; –0.3343) | –1.8897 | 578 | 502 | 15 |
| | | 2б | (1.0000; –0.3520) | –1.8877 | 557 | 449 | |
| | 0.1 | 1а | (0.3401; 1.0000) | –1.8904 | 8890 | 6064 | 16 |
| | | 1б | (0.3611; 1.0000) | –1.8820 | 10928 | 3551 | |
| | | 2а | (0.3409; 1.0000) | –1.8904 | 35511 | 34867 | 17 |
| | | 2б | (1.0000; –0.3399) | –1.8904 | 34965 | 34482 | |

Таблица 3: Результаты численных экспериментов, полученных Алгоритмом I

На рис. А.3–А.17 показаны сравнительные графики зависимости величины $h'$, существенно влияющей на объём вычислений, от номера $k$ итерации для различных вариантов проведённых тестов.

Из полученных результатов видно, что по количеству выполненных в процессе решения итераций среди вариантов реализации схемы 1 (обход в глубину) лучшим в большей части экспериментов является вариант 1а; среди вариантов реализации схемы 2 (обход в ширину) меньше итерационных шагов для большей части примеров требует вариант 2б.

Если же сравнивать схемы 1 и 2 между собой, то по количеству итераций, а также по машинному времени в подавляющем большинстве экспериментов выигрывает схема 1 (обход в глубину).

Результаты численных экспериментов, полученных алгоритмом II при различных ε и γ и β=0.99, представлены в табл. 4, в которой также указаны ссылки на соответствующие этим результатам рисунки. Через θ обозначено отношение числа итераций, где производилось деление на пять параллелепипедов, к числу итераций, где производилось деление пополам; прочерк означает, что метод не привел к решению за $10^5$ итераций.

| $i$ | ε | γ | $(x_\varepsilon, y_\varepsilon)$ | $f_i(x_\varepsilon, y_\varepsilon)$ | $N_{tot}$ | $N_{opt}$ | θ | Рис. А.№ |
|---|---|---|---|---|---|---|---|---|
| 1 | 0.5 | 0.01 | (0.00122, –0.00048) | –9.7119 | 427 | 24 | 2.298 | 18, 22 |
| | | 1 | (0.00122, –0.00048) | –9.7119 | 725 | 24 | 0 | |
| | 0.1 | 0.01 | $(-7.62\times10^{-6}, -6.10\times10^{-5})$ | –9.9415 | 1175 | 473 | 0.208 | 18, 23 |
| | | 1 | $(-7.62\times10^{-6}, -6.10\times10^{-5})$ | –9.9415 | 1337 | 473 | 0 | |
| 2 | 0.5 | 0.01 | (0.00122, 0.00033) | –12.4428 | 7613 | 6287 | 0.158 | 19, 24 |
| | | 1 | (0.00122, –0.00048) | –12.4301 | 9191 | 6418 | 0 | |
| | 0.1 | 0.01, 1 | – | – | – | – | – | – |
| 3 | 0.5 | 0.01 | (–9.48036, 9.48568) | –5.3339 | 33843 | 23855 | 0.062 | 20, 25 |
| | | 1 | (9.49218, –9.49218) | –5.3334 | 37975 | 34152 | 0 | |
| | 0.1 | 0.01, 1 | – | – | – | – | – | – |
| 4 | 0.5 | 0.01 | (–0.99778, 0.31519) | –1.8127 | 325 | 24 | ∞ | 21, 26 |
| | | 1 | (0.98437, –0.34375) | –1.7217 | 589 | 278 | 0 | |
| | 0.1 | 0.01 | (0.99724, –0.34830) | –1.8172 | 761 | 514 | 2.617 | 21, 27 |
| | | 1 | (0.99609, –0.34375) | –1.8045 | 1277 | 696 | 0 | |

Таблица 4: Результаты численных экспериментов, полученных Алгоритмом II

Из таблицы видно, что в некоторых случаях схема дробления, задаваемая значением γ, заметно влияла на эффективность метода. В частности, для функций $f_1$ и $f_4$ при ε = 0.5 и γ = 0.01 общее число параллелепипедов (оно совпадает с количеством вычислений значения целевой функции) примерно на 40 % меньше, чем при γ = 1 (производится только половинное деление). Подробная зависимость числа параллелепипедов $N_{tot}$ от параметра γ для функций $f_1$, …, $f_4$ представлена на рис. А.18–А.21.

На рис. А.22–А.27 показаны сравнительные графики зависимости величины $r_k$ от номера $k$ текущего параллелепипеда для значений γ и ε, представленных в табл. 4. Видно, что в каждом примере эти зависимости качественно одинаковы. Вместе с тем заметно преимущество дифференцированной схемы дробления, которое начинает проявляться не сразу, а только во второй трети итераций алгоритма или ближе к его концу.

# СПИСОК ЦИТИРОВАННОЙ ЛИТЕРАТУРЫ

# Приложение А

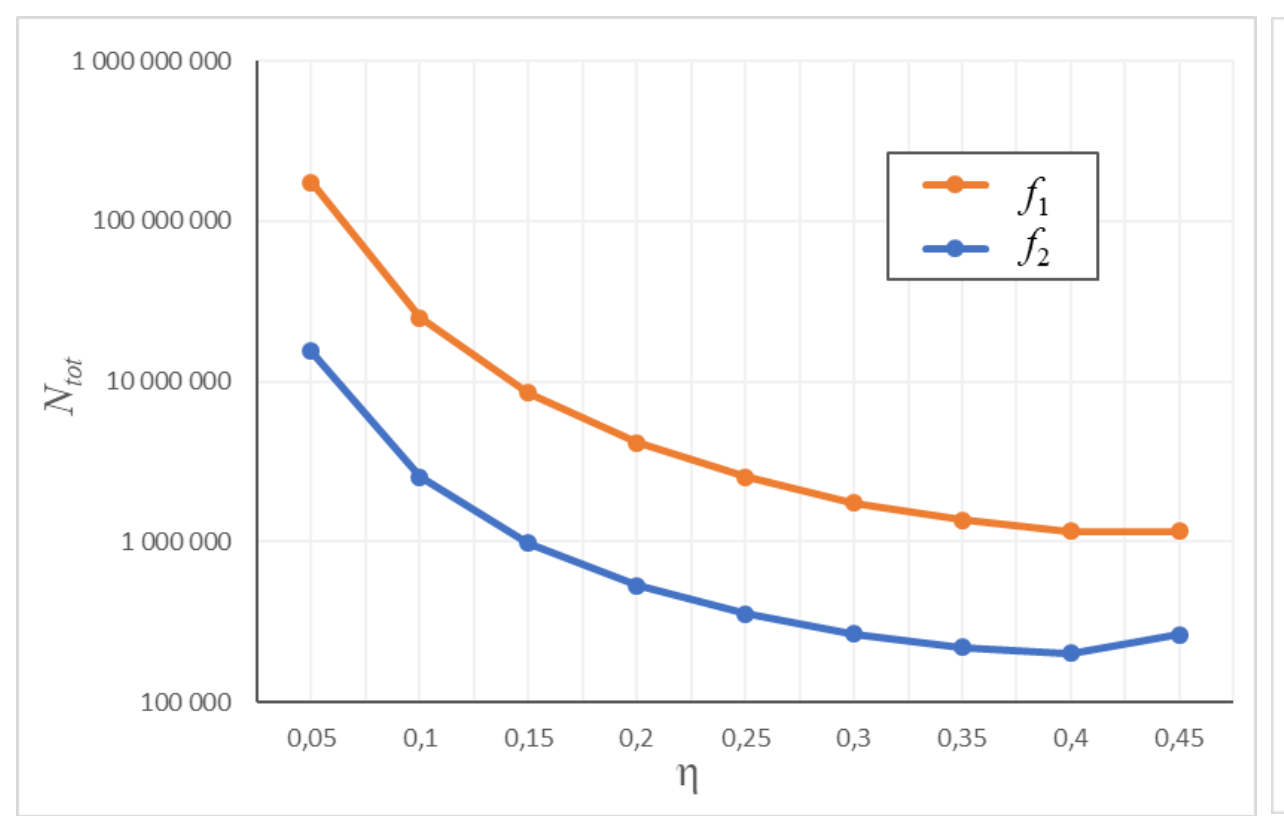

Рис. А.1. Зависимость $N_{tot}$ от η для $f_1$ и $f_2$ при ε=0.5

Рис. А.2. Зависимость $N_{tot}$ от η для $f_3$ и $f_4$ при ε=0.5

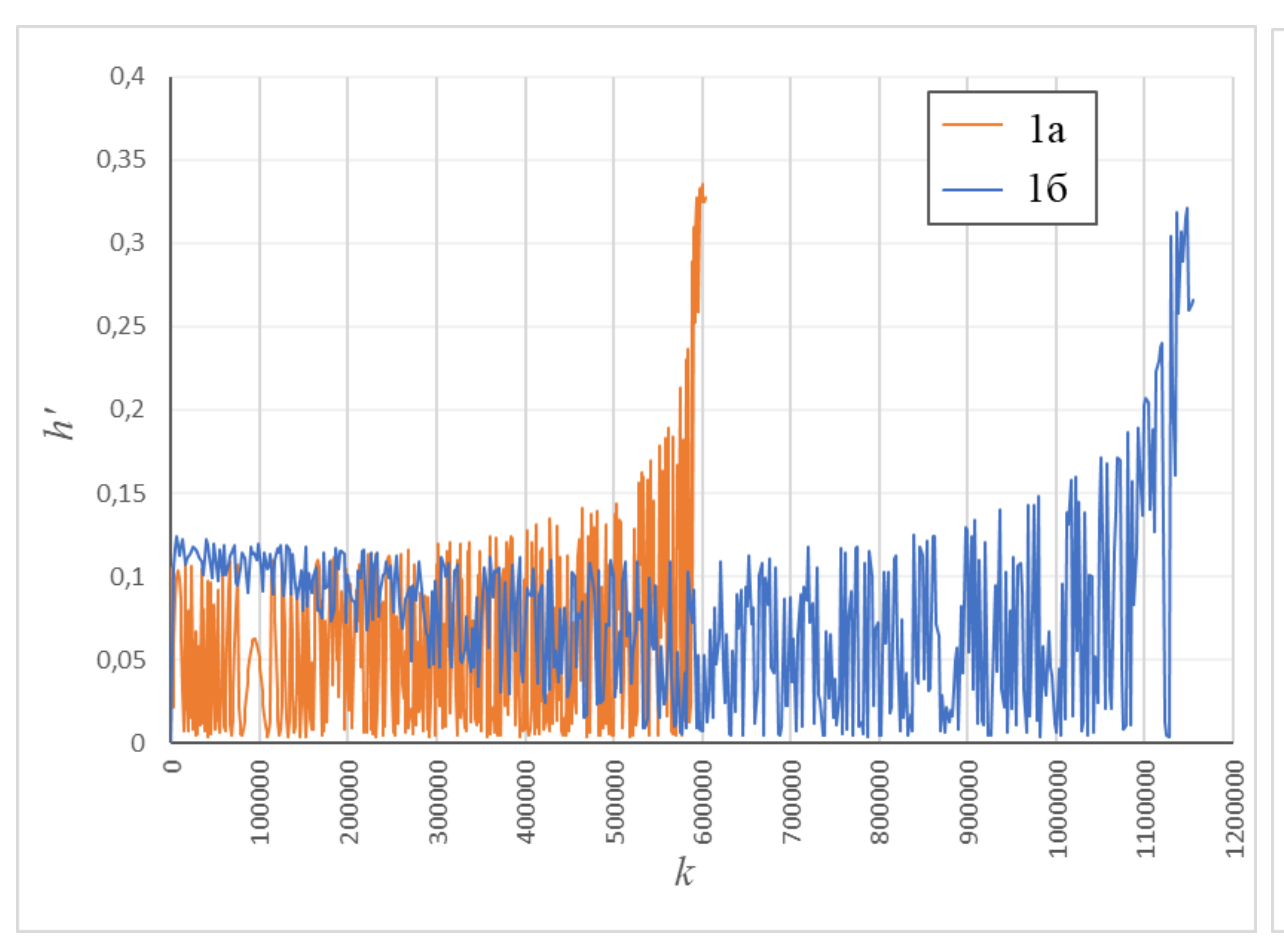

Рис. А.3. Зависимость $h'$ от $k$ для $f_1$ при ε=0.5 для схем 1а и 1б

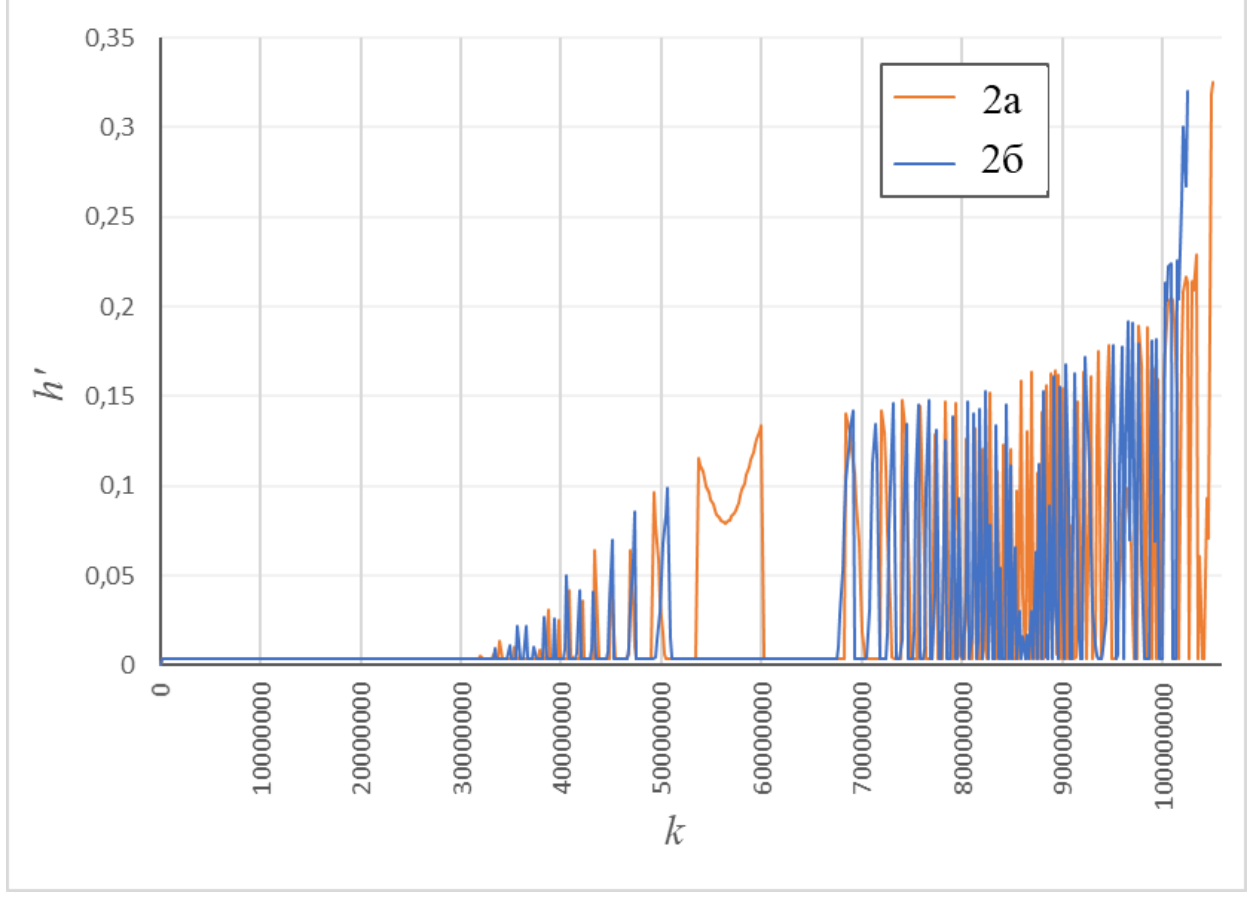

Рис. А.4. Зависимость $h'$ от $k$ для $f_1$ при ε=0.5 для схем 2а и 2б

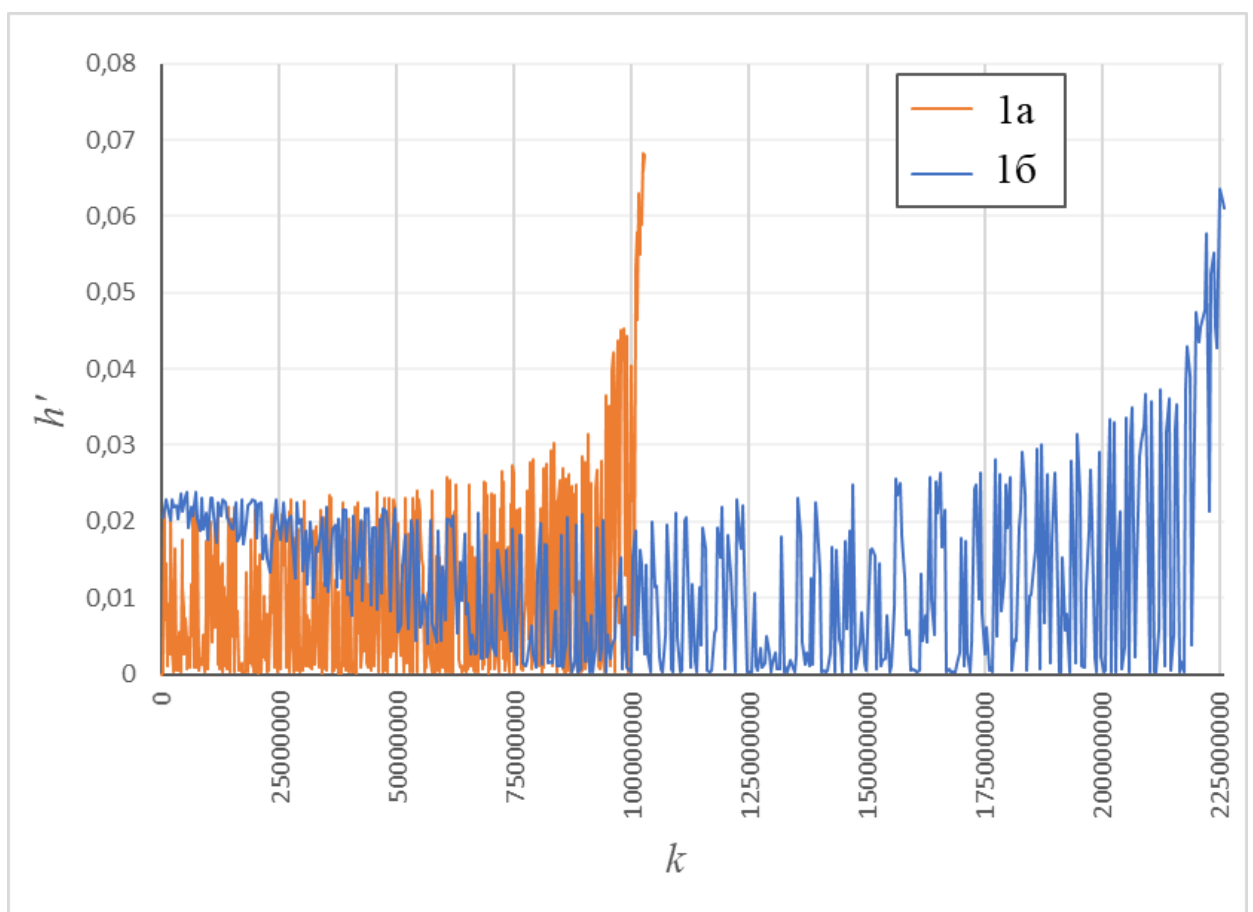

Рис. А.5. Зависимость $h'$ от $k$ для $f_1$ при ε=0.1 для схем 1а и 1б

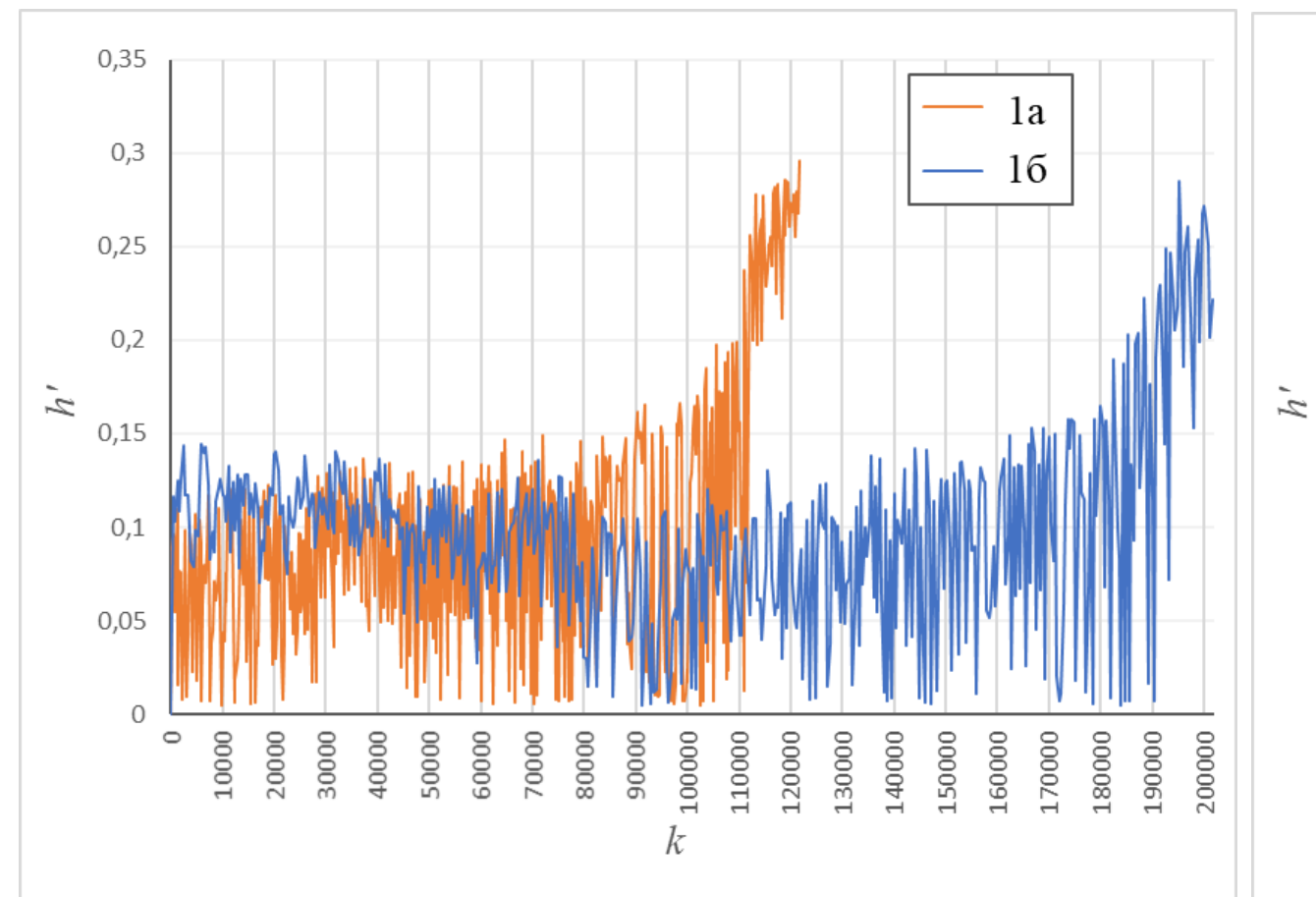

Рис. А.6. Зависимость $h'$ от $k$ для $f_2$
при ε=0.5 для схем 1а и 1б

Рис. А.7. Зависимость $h'$ от $k$ для $f_2$
при ε=0.5 для схем 2а и 2б

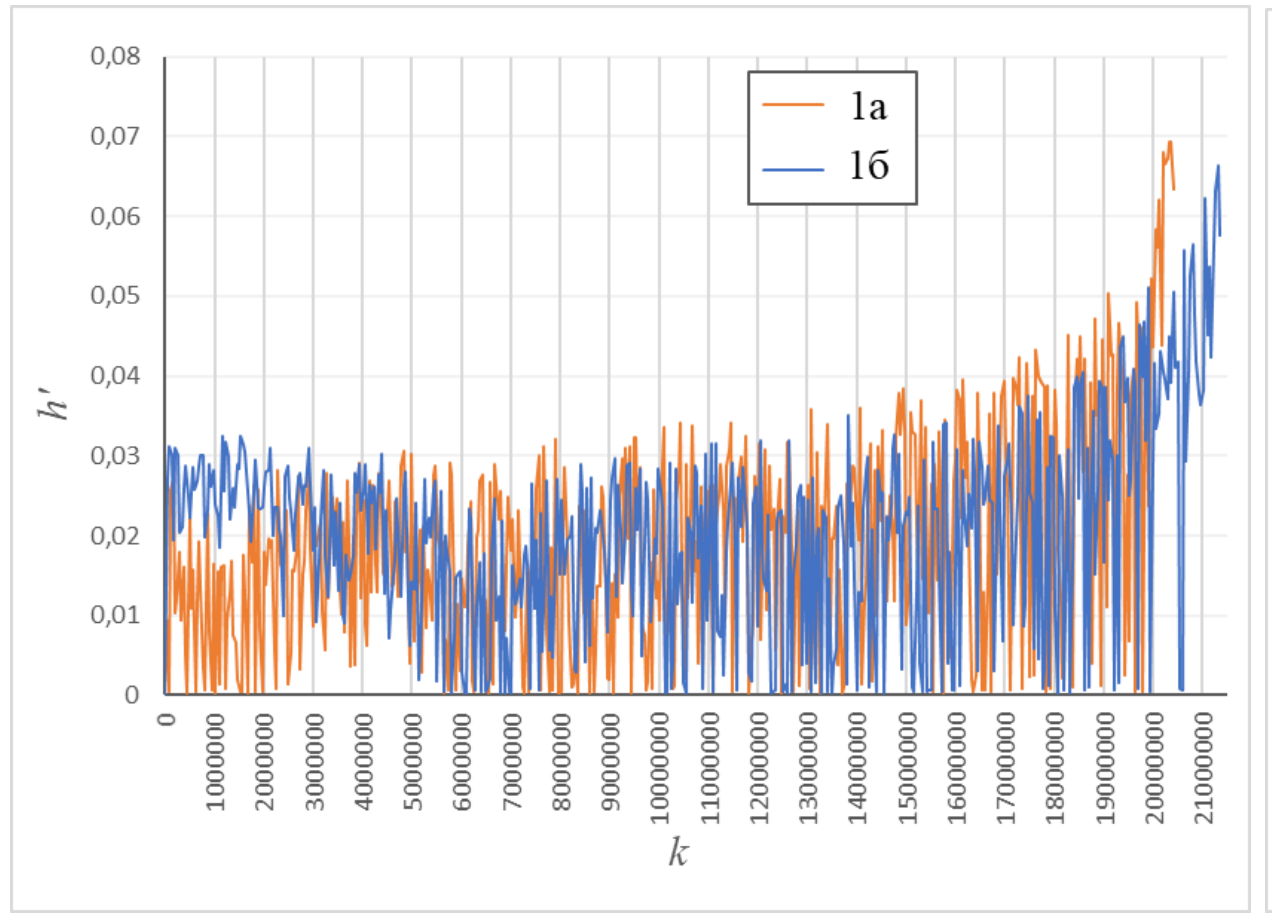

Рис. А.8. Зависимость $h'$ от $k$ для $f_2$
при ε=0.1 для схем 1а и 1б

Рис. А.9. Зависимость $h'$ от $k$ для $f_2$
при ε=0.1 для схем 2а и 2б

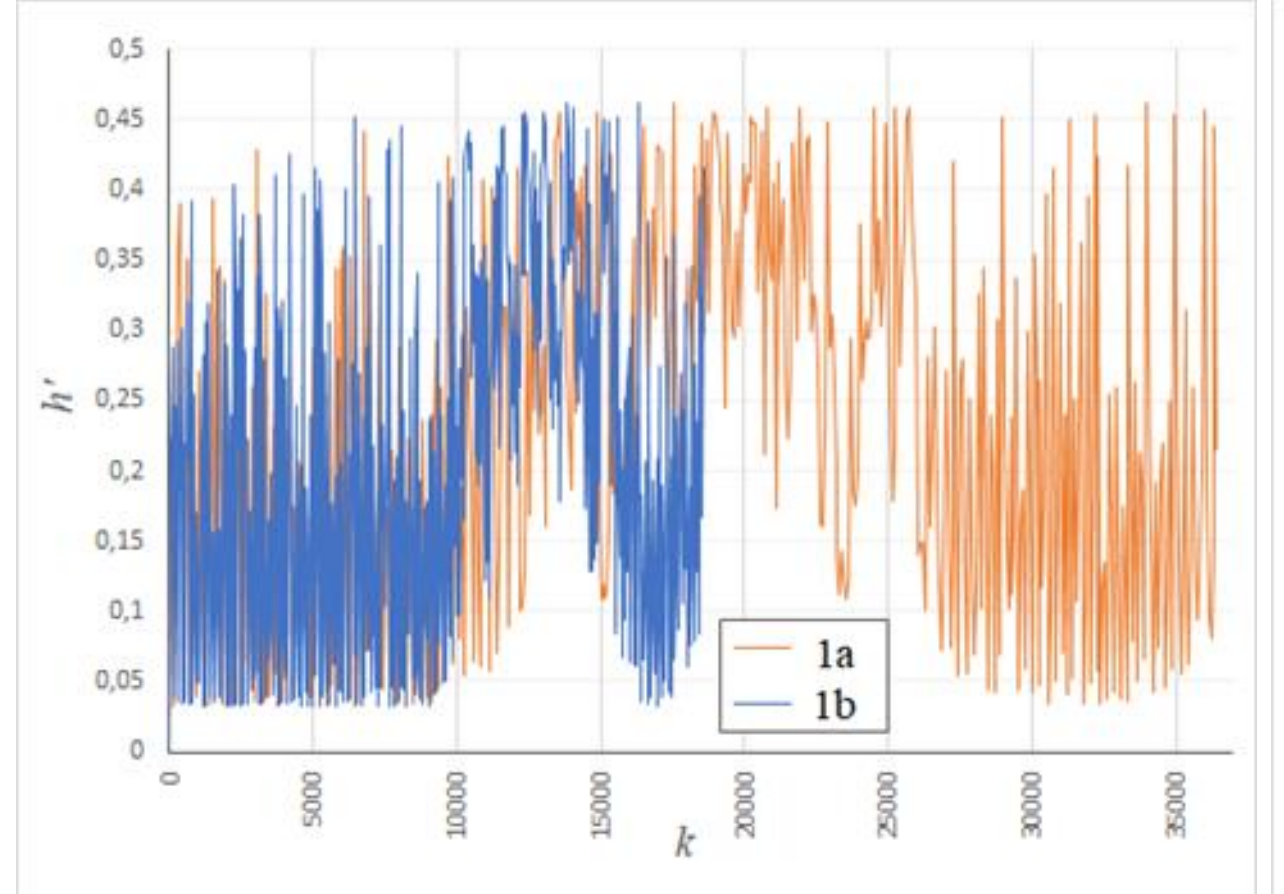

Рис. А.10. Зависимость $h'$ от $k$ для $f_3$
при ε=0.5 для схем 1а и 1б

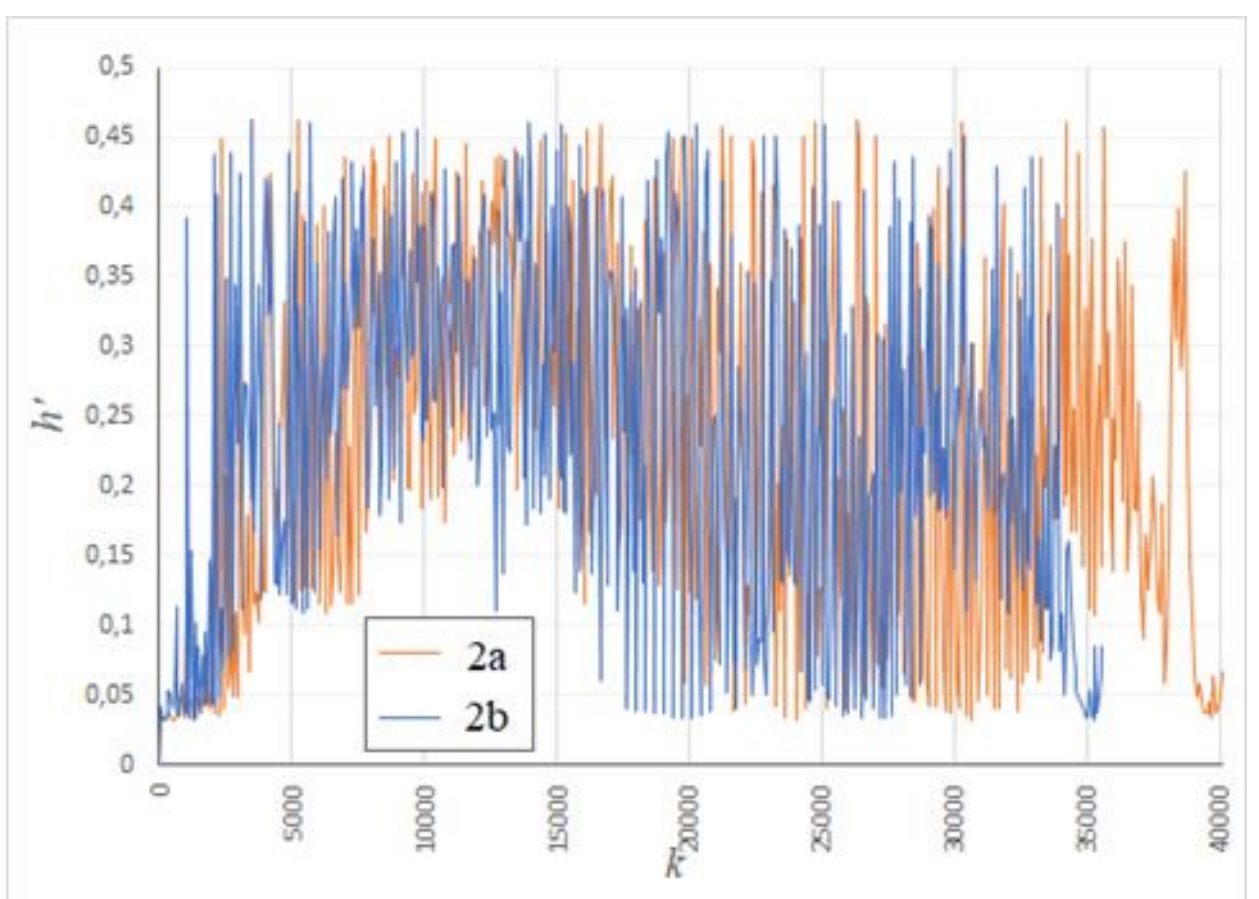

Рис. А.11. Зависимость $h'$ от $k$ для $f_3$
при ε=0.5 для схем 2а и 2б

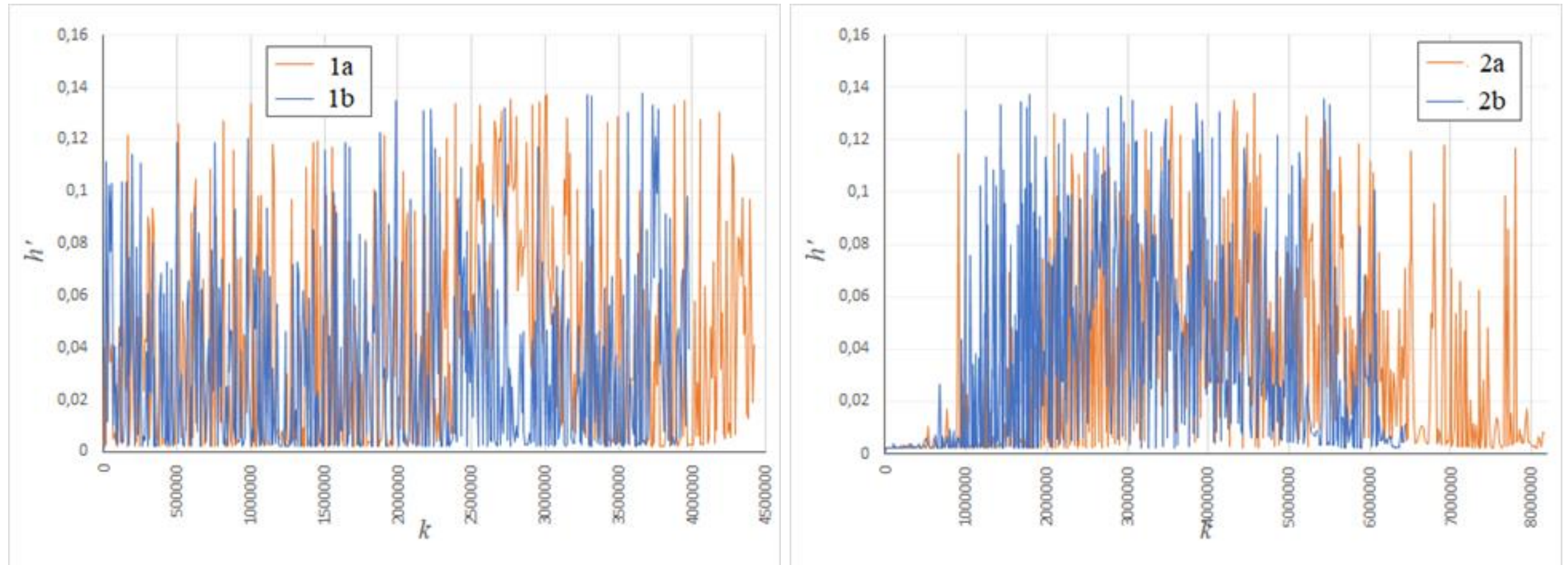

Рис. А.12. Зависимость $h'$ от $k$ для $f_3$ при ε=0.1 для схем 1а и 1б

Рис. А.13. Зависимость $h'$ от $k$ для $f_3$ при ε=0.1 для схем 2а и 2б

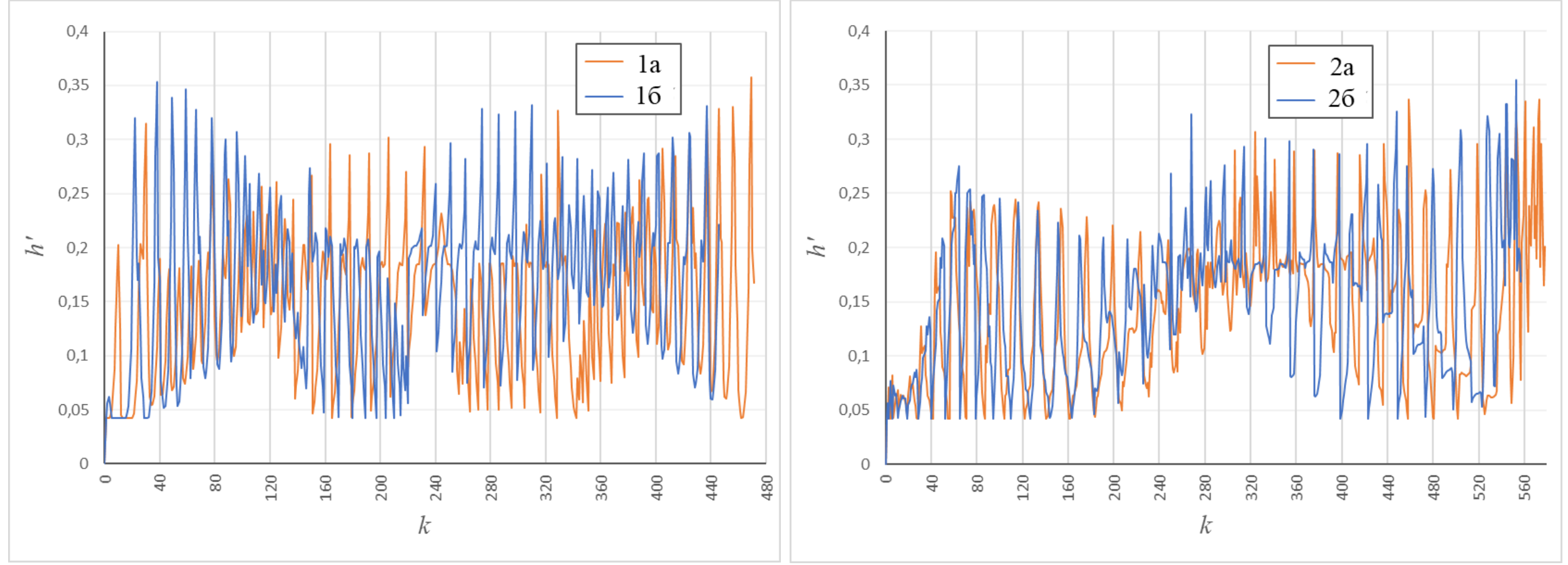

Рис. А.14. Зависимость $h'$ от $k$ для $f_4$ при ε=0.5 для схем 1а и 1б

Рис. А.15. Зависимость $h'$ от $k$ для $f_4$ при ε=0.5 для схем 2а и 2б

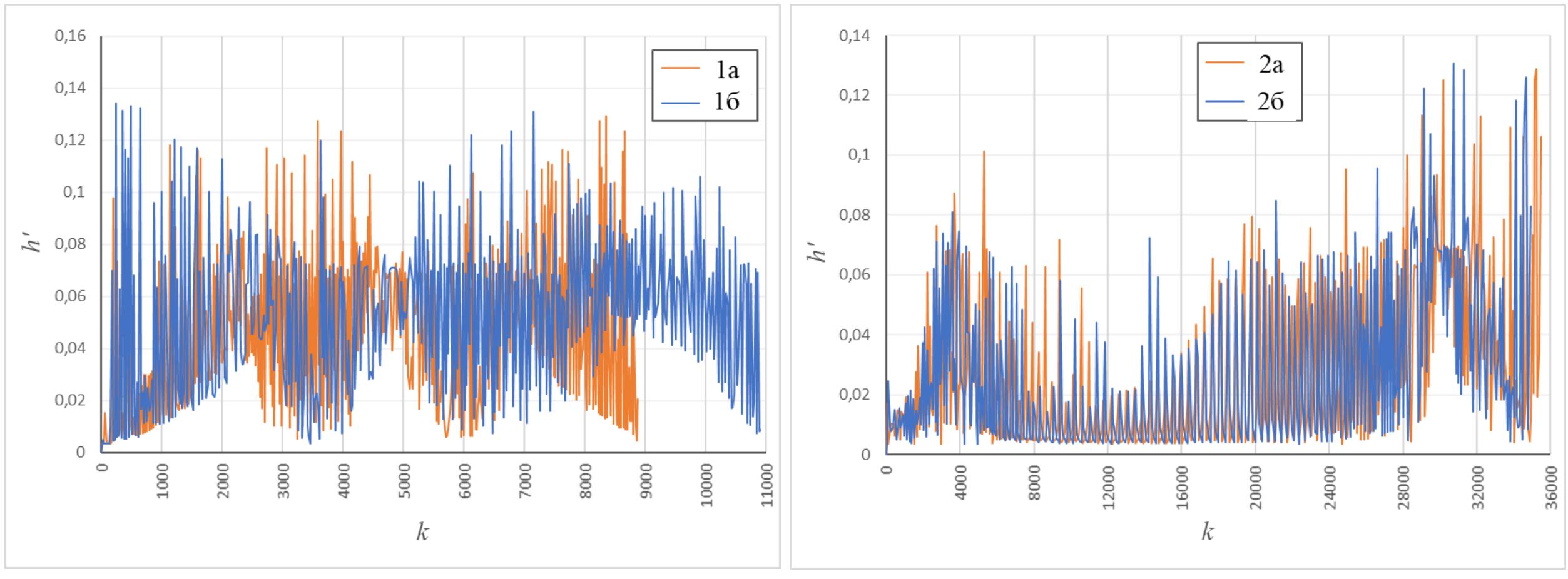

Рис. А.16. Зависимость $h'$ от $k$ для $f_4$ при ε=0.1 для схем 1а и 1б

Рис. А.17. Зависимость $h'$ от $k$ для $f_4$ при ε=0.1 для схем 2а и 2б

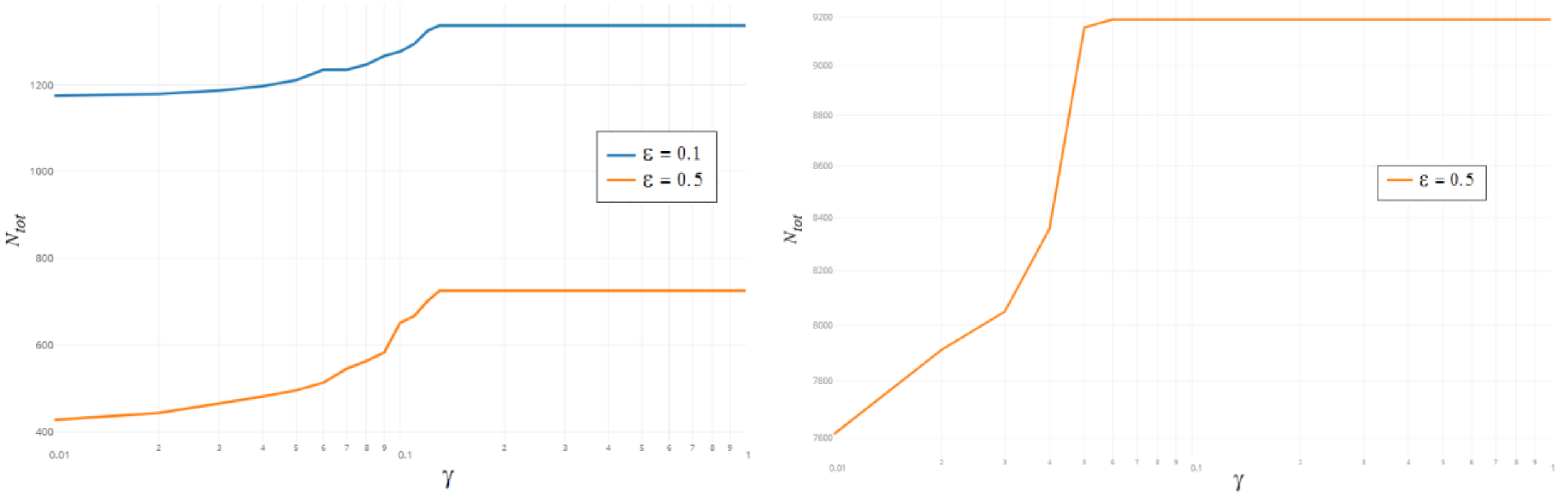

Рис. А.18 Зависимость $N_{tot}$ от γ для $f_1$

Рис. А.19. Зависимость $N_{tot}$ от γ для $f_2$

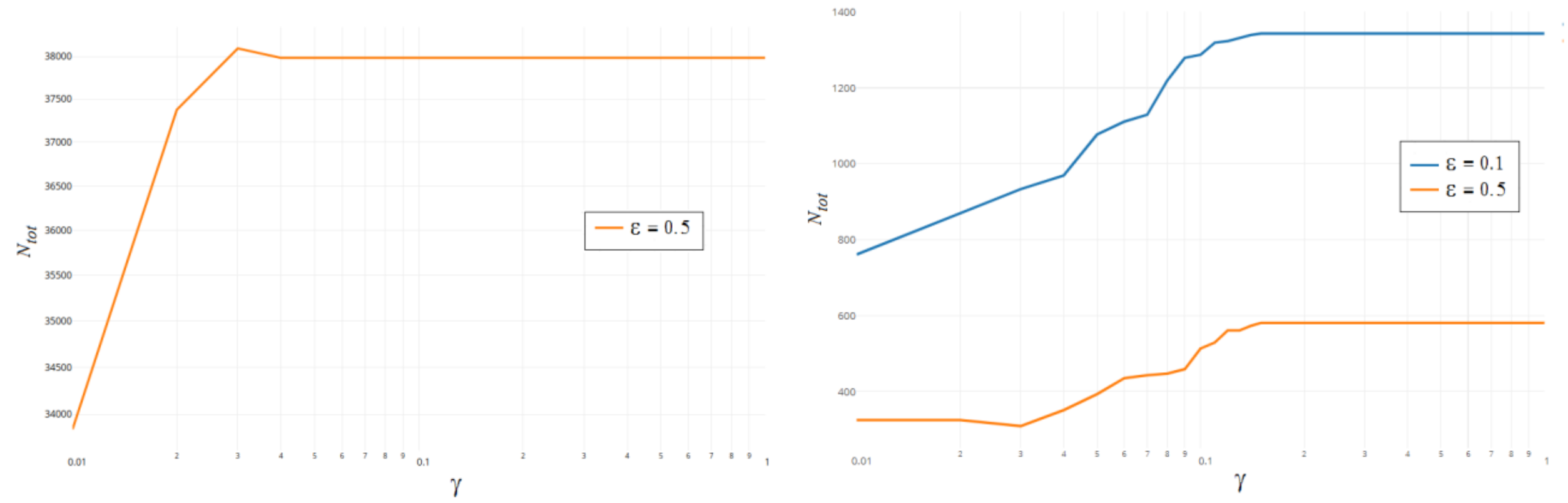

Рис. А.20. Зависимость $N_{tot}$ от γ для $f_3$

Рис. А.21. Зависимость $N_{tot}$ от γ для $f_4$

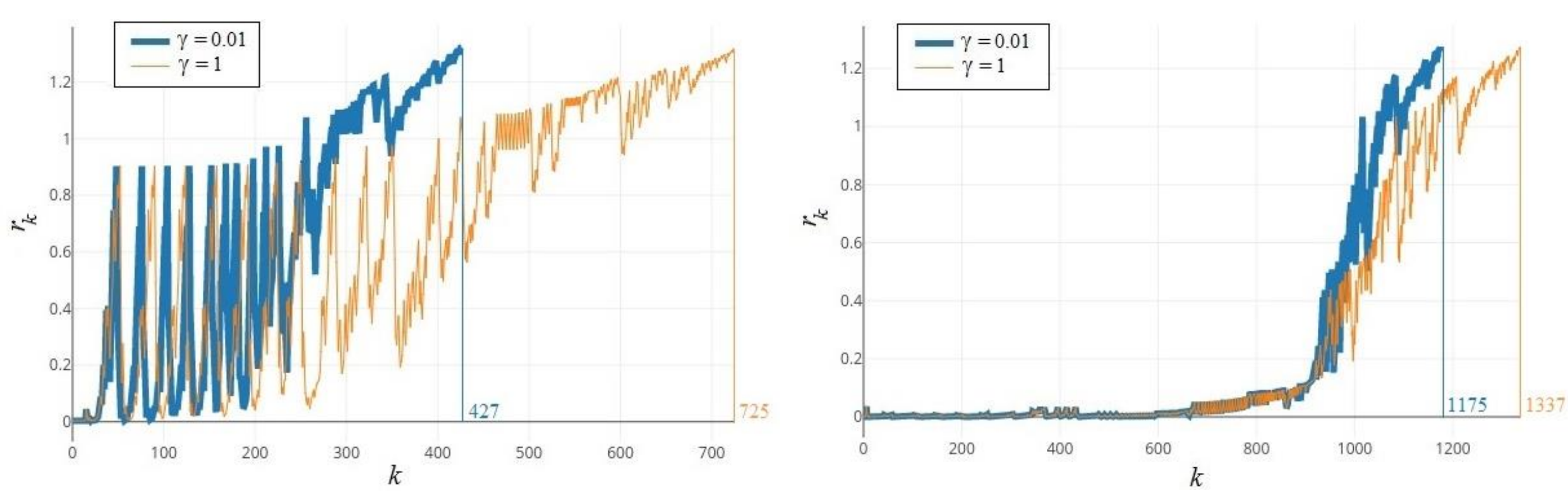

Рис. А.22. Зависимость $r_k$ от $k$ для $f_1$ при ε=0.5

Рис. А.23. Зависимость $r_k$ от $k$ для $f_1$ при ε=0.1

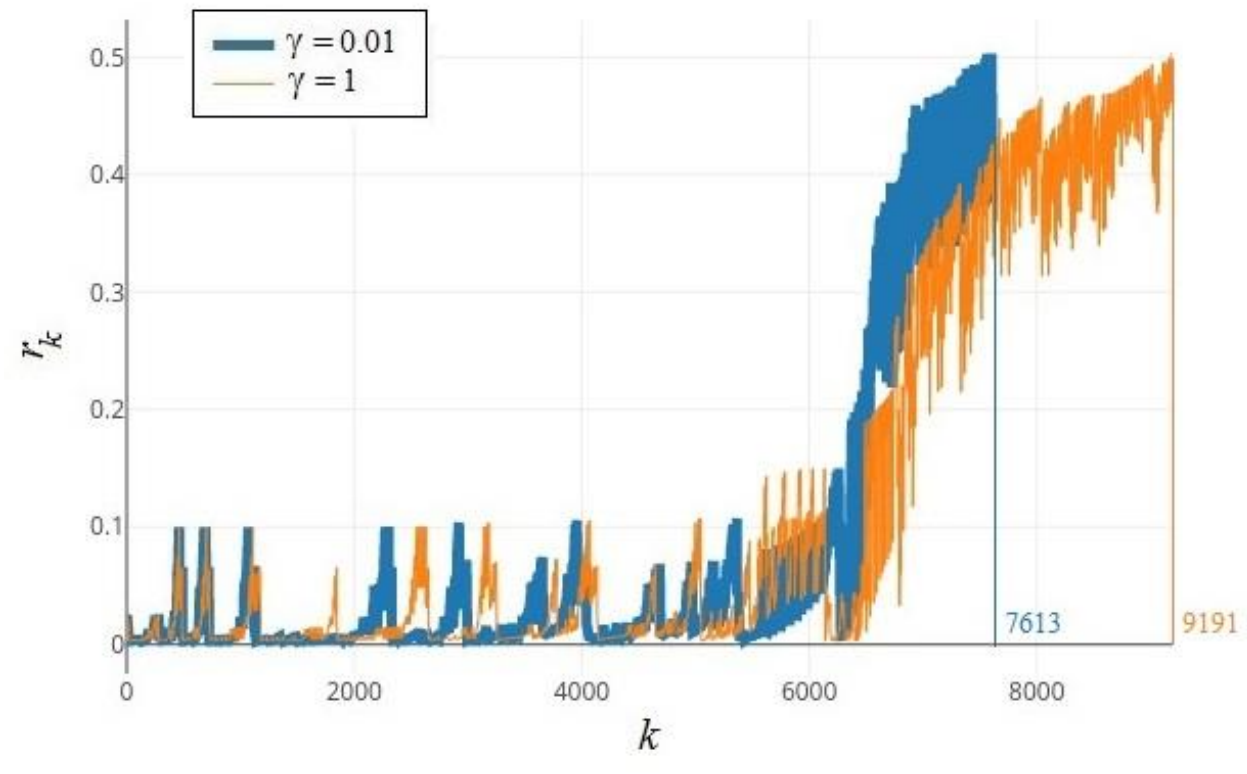

Рис. А.24. Зависимость $r_k$ от $k$ для $f_2$ при ε=0.5

Рис. А.25. Зависимость $r_k$ от $k$ для $f_3$ при ε=0.5

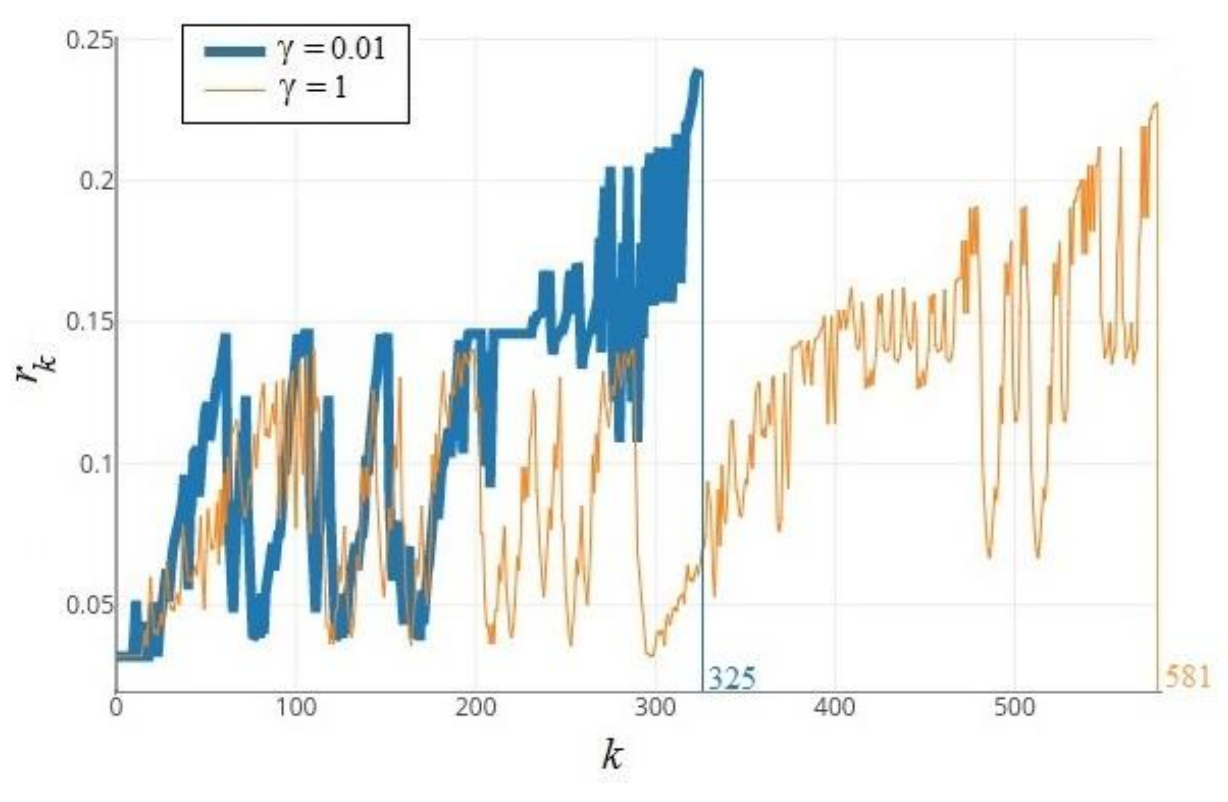

Рис. А.26. Зависимость $r_k$ от $k$ для $f_4$ при ε=0.5

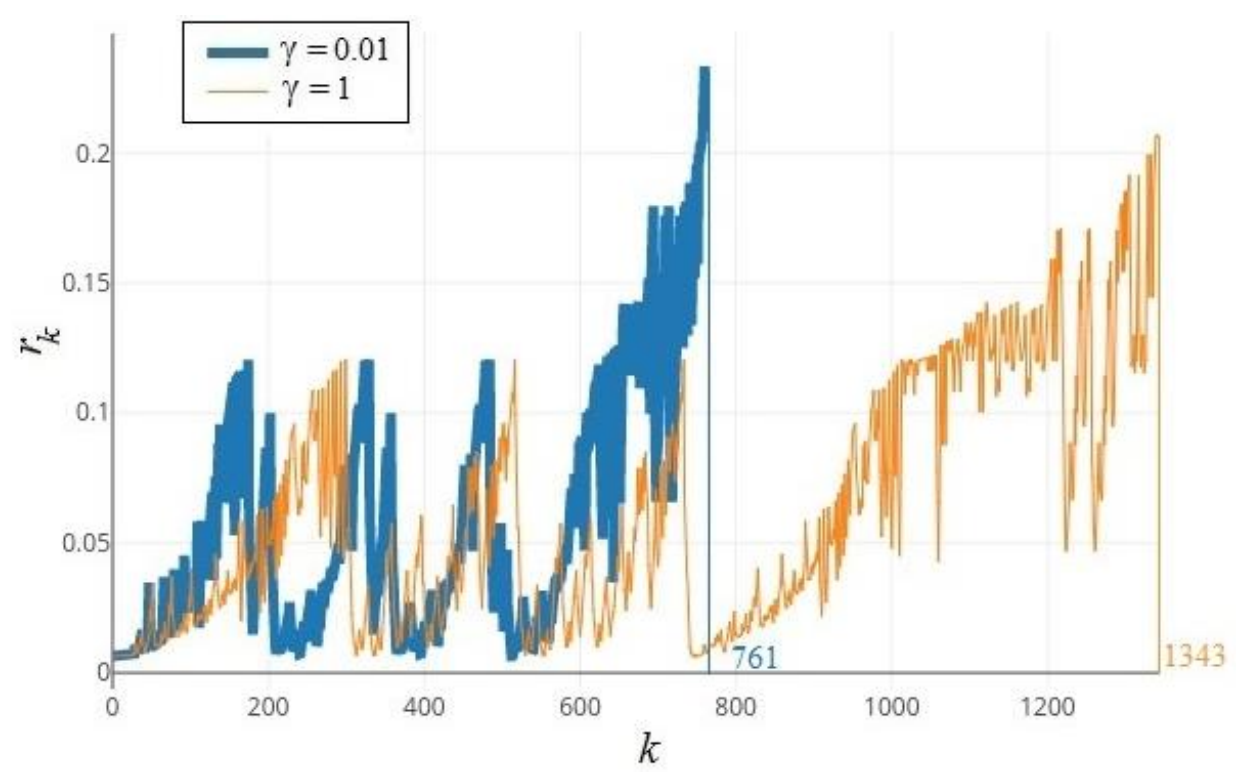

Рис. А.27. Зависимость $r_k$ от $k$ для $f_4$ при ε=0.1